\documentclass[12pt]{slides}
\usepackage{amsmath,amsfonts,amssymb,amsthm}
\begin{document}
\author{{\bf Denis~I.~Saveliev\/}}
\title{{\bf A Note on Singular Cardinals \\
in~Set~Theory~without~Choice\/}}
\date{2007 August 11, Beijing
\vfill{
\begin{tiny} 
Partially supported by grant~06-01-00608-a 
of Russian Foundation for Basic Research
\end{tiny}
}}
\maketitle

\theoremstyle{plain}
\newtheorem{thm}{Theorem}
\newtheorem{prb}{Problem}
\newtheorem{coro}{Corollary}
\newtheorem*{tm}{Theorem}
\newtheorem*{cor}{Corollary}
\newtheorem*{lm}{Lemma}
\newtheorem*{fct}{Fact}

\theoremstyle{definition}
\newtheorem*{df}{Definition}
\newtheorem*{rfr}{References}

\theoremstyle{remark}
\newtheorem*{rmk}{Remark}
\newtheorem*{exm}{Example}

\newcommand{\eqs}{ {\;=^{\!*}\,} }
\newcommand{\lesss}{ {\;<^{\!*}\,} }
\newcommand{\leqs}{ {\;\le^{\!*}\,} }
\newcommand{\wo}{\mathop {\mathrm {wo\,}}\nolimits }
\newcommand{\pwo}{\mathop {\mathrm {pwo\,}}\nolimits }
\newcommand{\wf}{\mathop {\mathrm {wf\,}}\nolimits }
\newcommand{\ewf}{\mathop {\mathrm {ewf\,}}\nolimits }
\newcommand{\cf}{ {\mathop{\mathrm {cf\,}}\nolimits} }
\newcommand{\tc}{ {\mathop{\mathrm {tc\,}}\nolimits} }
\newcommand{\dom}{ {\mathop{\mathrm {dom\,}}\nolimits} }
\newcommand{\ran}{ {\mathop{\mathrm {ran\,}}\nolimits} }
\newcommand{\fld}{ {\mathop{\mathrm {fld\,}}\nolimits} }
\newcommand{\otp}{ {\mathop{\mathrm {otp\,}}\nolimits} }
\newcommand{\lh}{ {\mathop{\mathrm {lh\,}}\nolimits} }
\newcommand{\gt}{\mathfrak}
\newcommand{\ext}{\mathrm{ext}}
\newcommand{\Z}{ {\mathrm Z} }
\newcommand{\ZF}{ {\mathrm {ZF}} }
\newcommand{\ZFA}{ {\mathrm {ZFA}} }
\newcommand{\ZFC}{ {\mathrm {ZFC}} }
\newcommand{\E}{ {\mathrm {AE}} }
\newcommand{\AR}{ {\mathrm {AR}} }
\newcommand{\WR}{ {\mathrm {WR}} }
\newcommand{\WO}{ {\mathrm {WO}} }
\newcommand{\LO}{ {\mathrm {LO}} }
\newcommand{\AF}{ {\mathrm {AF}} }
\newcommand{\AC}{ {\mathrm {AC}} }
\newcommand{\DC}{ {\mathrm {DC}} }
\newcommand{\AD}{ {\mathrm {AD}} }
\newcommand{\AInf}{ {\mathrm {AInf}} }
\newcommand{\AU}{ {\mathrm {AU}} }
\newcommand{\AP}{ {\mathrm {AP}} }
\newcommand{\CH}{ {\mathrm {CH}} }
\newcommand{\APr}{ {\mathrm {APr}} }
\newcommand{\ASp}{ {\mathrm {ASp}} }
\newcommand{\ARp}{ {\mathrm {ARp}} }
\newcommand{\AFA}{ {\mathrm {AFA}} }
\newcommand{\BAFA}{ {\mathrm {BAFA}} }
\newcommand{\FAFA}{ {\mathrm {FAFA}} }
\newcommand{\SAFA}{ {\mathrm {SAFA}} }
\newcommand{\NS}{ {\mathrm{NS}} }
\newcommand{\Cov}{ {\mathrm{Cov}} }

\newpage

In this talk,
I discuss 
{\it how singular\/} can cardinals be 
in absence of~AC, the axiom of choice. 
I~shall show that,
contrasting with known 
negative {\it consistency\/} results 
(of Gitik and others), 
certain positive results are {\it provable\/}. 
At the end,
I pose some problems.

\newpage

{\centerline{\Large{\bf Preliminaries\/}}}

\newpage

\noindent
\begin{df}
Given a set~$X$,
its {\it cardinal\/} number 
$$|X|$$
is the~class of all sets 
of~{\it the same size\/} that~$X$,
i.e.,
admitting a~one-to-one map onto~$X$.
\end{df}
Thus
$$|X|=|Y|$$
means 
``There is a~bijection of~$X$ onto~$Y$''.

\newpage

Cardinals of nonempty sets are proper classes;
so, we have a~little technical obstacle:
$$
\text{How quantify cardinals?}
$$
In some happy cases 
we can represent them by~sets:

If
$|X|$ is a~{\it well-ordered cardinal\/},
i.e.,
meets the~class of (von Neumann's) ordinals,
take the least such ordinal
(an~{\it initial\/} ordinal).

If
$|X|$  is a~{\it well-founded cardinal\/},
i.e.,
meets the~class of~well-founded sets,
take the~lower level of the~intersection
(so-called Scott's trick).

\newpage

What is in general?
The answer is 
$$
\text{No matter}
$$
because instead of cardinals,
we can say about sets and bijections.

Thus
$$
\varphi(|X|,|Y|,\ldots)
$$
means
$
\varphi(X',Y',\ldots)
$
whenever $|X|=|X'|$, $|Y|=|Y'|$,~$\ldots$

\newpage

Notations:

The German letters
$$\gt l,\gt m,\gt n,\ldots$$\!
denote arbitrary cardinals.
The Greek letters
$$\lambda,\mu,\nu,\ldots$$\!
denote well-ordered ones
(i.e., initial ordinals),
while
the Greek letters
$$\alpha,\beta,\gamma,\ldots$$\!
denote arbitrary ordinals.

\newpage

Two basic relations on~cardinals
(dual in a~sense):
$$
|X|\le|Y|
$$
means 
``$X$ is empty or there is an injection of~$X$ into~$Y$'',
and
$$
|X|\leqs|Y|
$$
means
``$X$ is empty or there is a~surjection of~$Y$ onto~$X$''.

Equivalently,

$|X|\le|Y|$ means
``There is a~subset of~$Y$ of size~$|X|$'',

$|X|\leqs|Y|$ means
``$X$ is empty or there is a~partition of~$Y$ into $|X|$~pieces''.

\newpage

Clearly:

(i)
Both $\le$ and $\!\leqs\!$ are reflexive and transitive.

(ii)
$\le$ is antisymmetric
(Dedekind; Bernstein),
$\leqs\!$ is not necessarily.

(iii)
$\le$ is stronger than~$\!\leqs\!$.
Both relations coincide on well-ordered cardinals.

\newpage

Two important functions on cardinals
(Hartogs and Lindenbaum resp.):
$$
\aleph(\gt n)=
\{\alpha:|\alpha|\le\gt n\},
$$
$$
\aleph^*(\gt n)=
\{\alpha:|\alpha|\leqs\gt n\}.
$$

Equivalently,

$\aleph(\gt n)$
is the least~$\alpha$
such that
on a~set of size~$\gt n$
there is no well-ordering
of length~$\alpha$,

$\aleph^*(\gt n)$ 
is the least~$\alpha$
such that
on a~set of size~$\gt n$
there is no pre-well-ordering
of length~$\alpha$.

Customarily,
$\nu^+$ denotes $\aleph(\nu)$ 
for $\nu$ well-ordered.

\newpage

Clearly:

(i)
$\aleph(\gt n)$ 
and 
$\aleph^*(\gt n)$
are well-ordered cardinals.

(ii)
$\aleph(\gt n)\not\leq\gt n$ 
and 
$\aleph^*(\gt n)\not\!\!\!\!\leqs\gt n.$

It follows
$\nu<\nu^+$
and so
$$
\aleph_0<\aleph_1<\ldots<
\aleph_\omega<\ldots<\aleph_{\omega_1}<\ldots
$$
(where $\aleph_\alpha$ is 
$\alpha$th~iteration of~$\aleph$
starting from~$\aleph_0$).

(iii)
$
\aleph(\gt n)\le
\aleph^*(\gt n),
$
and both operations coincide
on well-ordered cardinals.
On other cardinals,
the gap can be very large:

\noindent
\begin{exm}
Assume~$\AD$.
Then
$\aleph(2^{\aleph_0})=\aleph_1$
while
$\aleph^*(2^{\aleph_0})$
is a~very large cardinal
(customarily denoted~$\Theta$).
\end{exm}

\newpage

{\centerline{\Large{\bf 
Results on Singularity\/}}}

\newpage

Notations:
$$\Cov(\gt l,\gt m,\gt n)$$
means 
``A~set of size~$\gt n$ can be covered by
$\gt m$ sets of size~$\gt l$\,''.

$\Cov(\!<\!\gt l,\gt m,\gt n)$
and
$\Cov(\gt L,\gt m,\gt n)$
(where $\gt L$ is a~class of cardinals)
have the appropriate meanings.

\noindent
\begin{df}
$\,$
A~cardinal $\gt n$ is {\it singular\/} 
iff\\
$\Cov(\!<\!\gt n,<\!\gt n,\gt n)$,
and {\it regular\/} otherwise.
\end{df}

\newpage

What is under~AC?

\noindent
\begin{fct}
$\,$
Assume~$\AC$.
Then
$\Cov(\gt l,\gt m,\gt n)$
implies
$\gt n\le\gt l\cdot\gt m$.
\end{fct}

\noindent
\begin{cor}
$\,$
Assume~$\AC$.
Then all the successor alephs are regular.
\end{cor}

Thus 
$\neg\Cov(\lambda,\lambda,\lambda^+)$
\,for all $\lambda\ge\aleph_0$.

\newpage

What happens without~AC?

\noindent
\begin{tm}[Feferman L\'evy]
$\,$
$\aleph_1$~can be singular.
\end{tm}

Thus
$\Cov(\aleph_0,\aleph_0,\aleph_1)$
is consistent.

Moreover,
under a~large cardinal hypothesis,
so can be all uncountable alephs:

\noindent
\begin{tm}[Gitik]
$\,$
All uncountable alephs can be singular.
\end{tm}

Clearly, 
then 
$\Cov(\!<\!\lambda,\aleph_0,\lambda)$
for all $\lambda\ge\aleph_0$.

\newpage

\begin{rmk}
What is the consistency strength?

Without successive singular alephs:
\rule{0mm}{0mm}
\\
The same as of~ZFC.

With $\lambda,\lambda^+$ both singular:
\rule{0mm}{0mm}
\\
Between 1~Woodin cardinal 
(Schindler improving Mitchell)
and $\omega$~Woodin cardinals 
(Martin Steel Woodin).

So, in general case:
\rule{0mm}{0mm}
\\
A~proper class of Woodins.
\end{rmk}

\newpage

Specker's problem:

Is
$\Cov(\aleph_\alpha,\aleph_0,2^{\aleph_\alpha})$
consistent
for all $\alpha$ simultaneously?

Partial answer:

\noindent
\begin{tm}[Apter~Gitik]
$\,$
Let
$A\subseteq Ord$ consist
either

(i) of all successor ordinals; or

(ii) of all limit ordinals and all successor ordinals of form 
$\alpha=3n, 3n+1, \gamma+3n$, or $\gamma+3n+2$,
where $\gamma$ is a limit ordinal.

Then
$$
(\forall\alpha\in A)\,
\Cov(\aleph_\alpha,\aleph_0,2^{\aleph_\alpha})
$$ 
is consistent (modulo large cardinals).
\end{tm}

(Really, their technique gives slightly more.)\\
In general, the problem remains open.

\newpage

Question:
\,{\it How singular\/} can cardinals be without~AC?
in the following sense:
How small are $\gt l\le\gt n$ and $\gt m\le\gt n$ satisfying

(i)
$\Cov(\!<\!\gt l,<\!\gt n,\gt n)$?

(ii)
$\Cov(\!<\!\gt n,<\!\gt m,\gt n)$?

(iii)
$\Cov(\!<\!\gt l,<\!\gt m,\gt n)$?

On~(iii):
\rule{0mm}{0mm}
\\
Specker's problem is a~partial case.

On~(ii):
\rule{0mm}{0mm}
\\
The answer is
$$
\text{As small as possible}
$$
since
Gitik's model satisfies
$\Cov(\!<\!\gt n,\aleph_0,\gt n)$
for all (not only well-ordered)~$\gt n$.

\newpage

On~(i):
\rule{0mm}{0mm}
\\
For well-ordered~$\gt n$,
the answer is
$$
\gt l<\gt n\text{ is impossible.}
$$

\noindent
\begin{thm}
$\,$
$\Cov(\!<\!\lambda,\gt m,\nu)$
\,implies\,
$
\nu\leqs\lambda\cdot\gt m,
$
\,and so
$$
\nu^+\le\aleph^*(\lambda\cdot\gt m).
$$
\end{thm}

\noindent
\begin{cor}
$\,$
$\neg\Cov(\!<\!\lambda,\lambda,\lambda^+)$
\,for all $\lambda\ge\aleph_0$.
\end{cor}

Since
$\Cov(\lambda,\lambda,\lambda^+)$
is consistent,
the result is exact.

\noindent
\begin{rmk}
$\neg\Cov(\aleph_0,\aleph_0,\aleph_2)$
is an old result of Jech.
(I~am indebted to Prof.~Blass
who informed~me.)
By Corollary, really
$\neg\Cov(\aleph_0,\aleph_1,\aleph_2)$.
\end{rmk}

\newpage

\par
Next question:
Let $\Cov(\gt l,\gt m,\gt n)$,
is $\gt n$ estimated via $\gt l$ and~$\gt m$?
(when $\gt n$ is not well-ordered).
Without Foundation,
the answer is
$$
\text{No}
$$
Even in the simplest case
$\gt l=2$ and $\gt m=\aleph_0$
such an estimation of~$\gt n$
is not provable:

\noindent
\begin{thm}
$\,$
It is consistent that
for any $\gt p$
there exists
$\gt n\nleq\gt p$
such that
\,$\Cov(2,\aleph_0,\gt n)$.
\end{thm}

The proof uses 
a~generalization of
permutation model technique
to the case of a~proper class of atoms.
We use non-well-founded sets
instead of atoms.

\newpage

On the other hand,
$\aleph(\gt n)$ and $\aleph^*(\gt n)$
are estimated via
$\aleph(\gt l)$,
$\aleph^*(\gt l)$,
and
$\aleph^*(\gt m)$:

\noindent
\begin{thm}
\rule{0mm}{0mm}
\\
$\Cov(\gt L,\gt m,\gt n)$
\,implies\,
$$
\aleph(\gt n)\le
\aleph^*(\sup_{\gt l\in\gt L}\,\aleph(\gt l)\cdot\gt m)
$$
and
$$
\aleph^*(\gt n)\le
\aleph^*(\sup_{\gt l\in\gt L}\,\aleph^*(\gt l)\cdot\gt m).
$$
\end{thm}

\noindent
\begin{coro}
\rule{0mm}{0mm}
\\
$
\neg\Cov(\!<\!\lambda,\lambda,2^\lambda)
$
\,and\,
$
\neg\Cov(\gt n,2^{\gt n^2},2^{2^{\gt n^2\cdot 2}}).
$
\end{coro}

In particular:

\noindent
$\neg\Cov(\lambda,2^\lambda,2^{2^\lambda})$
and
$\neg\Cov(\beth_\alpha,\beth_{\alpha+1},\beth_{\alpha+2}).$

Since
$\Cov(\gt n,\gt n,2^\gt n)$
is consistent,
the result is near optimal.

\newpage

Another corollary is that
Specker's request,
even in a~weaker form, gives 
the least possible evaluation of
$\aleph^*(2^\lambda)$
(which is~$\lambda^{++}$):

\noindent
\begin{coro}
$\,$
$\Cov(\lambda,\lambda^+,2^\lambda)$
\,implies\,
$$
\aleph^*(2^\lambda)=
\aleph(2^\lambda)=
\lambda^{++}.
$$        
\end{coro}

So,
if there exists a~model 
which gives the positive answer to Specker's problem,
then in it, 
all the cardinals~$\aleph^*(2^\lambda)$
have the least possible values.

\newpage

As the last corollary,
we provide a~``pathology'' 
when a~set admits 
{\it neither\/} 
well-ordered covering (of arbitrary size)
by sets of smaller size,
{\it nor\/} 
covering of smaller size
by well-orderable sets (of arbitrary size).
Moreover,
it can be the {\it real line\/}:

\noindent
\begin{coro}
$\,$
Assume
$\CH$ holds and $\Theta$~is limit.
(E.g., assume $\AD$.)
Then
for any well-ordered~$\lambda$
$$
\neg\Cov(\!<\!2^{\aleph_0},\lambda,2^{\aleph_0})
\text{\, and \,}
\neg\Cov(\lambda,\!<\!2^{\aleph_0},2^{\aleph_0}).
$$
\end{coro}

(Here CH means
``There is no~$\gt m$ such that\\
$\aleph_0<\gt m<2^{\aleph_0}$\,''.)

\newpage

{\centerline{\Large{\bf 
Problems\/}}}

\newpage

\noindent
\begin{prb}
$\,$
Is\, 
$\neg\Cov(\gt n,2^\gt n,2^{2^\gt n})$
true
for all~$\gt n$?
\end{prb}

That holds if $\gt n=\gt n^2$
(by Corollary~1 of Theorem~3).

\noindent
\begin{prb}
$\,$
Is\,
$\neg\Cov(\!<\!\beth_\alpha,\beth_\alpha,\beth_{\alpha+1})$
true
for all~$\alpha$?
\end{prb}

That near holds if $\alpha$~is successor
(again by Corollary~1 of Theorem~3).

\noindent
\begin{prb}
$\,$
Is\,
$\Cov(\gt n,\aleph_0,2^{\gt n^2})$
consistent
for all~$\gt n$ simultaneously?
\end{prb}

This sharps Specker's problem of course.

\noindent
\begin{prb}
$\,$
Can Theorem~2 be proved
assuming Foundation?
More generally,
expand the Transfer Theorem (Jech Sohor)
to the case of a~proper class of atoms.
\end{prb}

\newpage

\noindent
\begin{prb}
$\,$
Is it true that
on successor alephs
the cofinality can behave anyhow,
in the following sense:
Let
$F$ be any function
such that
$$
F:SuccOrd\to SuccOrd\cup\{0\}
$$
and $F$ satisfies
$$
\begin{array}{rll}
\text{{\rm (i)\/}}&
F(\alpha)\le\alpha\;\text{ and}
\\
\text{{\rm (ii)\/}}&
F(F(\alpha))=F(\alpha)
\end{array}
$$
for all successor~$\alpha$.
Is it consistent
$$
\cf\aleph_\alpha=\aleph_{F(\alpha)}
$$
for all successor~$\alpha$?
\end{prb}

Perhaps
if~$F$ makes no successive cardinals singular,
it is rather easy;
otherwise very hard.

\newpage

{\centerline{\Large{\bf 
References\/}}}

\newpage

[1]~Arthur~W.~Apter
and
Moti~Gitik.
{\it Some results on Specker's problem\/}.
Pacific Journal of~Mathematics,
134,~2~(1988),
227--249.

[2]~Solomon~Feferman
and
Azriel~L\'evy.
{\it Independences results
in set theory by Cohen's method,~II.\/}
Notices of the American Mathematical Society,
10~(1963),
593.
Abstract.

[3]~Moti~Gitik.
{\it All uncountable cardinals can be singular\/}.
Israel Journal of~Mathematics,
35,~1--2~(1980),
61--88.

[4]~Moti~Gitik.
{\it Regular cardinals in models of~$\ZF$\/}.
Transactions of the~American Mathematical Society,
290,~1~(1985),
41--68.

[5]~Donald~A.~Martin
and
John~R.~Steel.
{\it Projective determinacy\/}.
Proceedings of the National Academy of Sciences of~U.S.A.,
85,~18 (1988),
6582--6586.

[6]~Donald~A.~Martin
and
John~R.~Steel.
{\it A~proof of projective determinacy\/}.
Journal of the American Mathematical Society,
2,~1~(1989),
71--125.

[7]~Ralf~Dieter~Schindler.
{\it Successive weakly compact or singular cardinals\/}.
Journal of~Symbolic Logic,
64~(1999),
139--146.

[8]~Ernst~P.~Specker.
{\it Zur Axiomatik der Mengenlehre (Fundierungs- und Auswahlaxiom)\/}.
\\Zeitschrift f\"ur~Mathematische Logik und Grundlagen der Mathematik,
3,~3~(1957),
173--210.

[9]~W.~Hugh~Woodin.
{\it Supercompact cardinals, sets of reals, and
weakly homogeneous trees\/}.
Proceedings of the National Academy of Sciences of~U.S.A.,
85,~18 (1988),
6587--6591.

\end{document}